\documentclass{amsart}
\usepackage{amssymb,amsmath}
\usepackage{amsthm}
\usepackage{euscript}
\usepackage{arydshln}
\input{epsf}

\newcounter{sec}

\newcounter{punct}[sec]

\def\punct{\refstepcounter{punct}{\arabic{section}.\arabic{punct}.  }}

\newtheorem{theorem}{Theorem}[sec]
\newtheorem{proposition}[theorem]{Proposition}

\newtheorem{lemma}[theorem]{Lemma}

\newtheorem{corollary}[theorem]{Corollary}

\def\COUNTERS{\addtocounter{sec}{1}
              \setcounter{punct}{0}
          \setcounter{equation}{0}
          \setcounter{theorem}{0}
         }
          
          \def\sm{\smallskip}

\newcommand{\gr}{\mathop {\mathrm {gr}}\nolimits}
\newcommand{\rk}{\mathop {\mathrm {rk}}\nolimits}

\def\0{\mathbf 0}

\def\ov{\overline}
\def\wh{\widehat}

\renewcommand{\rk}{\mathop {\mathrm {rk}}\nolimits}

 \def\C {{\mathbb C }}

\def\P{\mathbb P}

\def\M{\text{\sf M}}

\def\cP{\EuScript P}

\def\cA{\EuScript A}
\def\cO{\EuScript O}
\def\cT{\EuScript T}

\def\bbA{\mathbb A}

\def\kappa{\varkappa}
\def\epsilon{\varepsilon}
\def\phi{\varphi}
\def\le{\leqslant}
\def\ge{\geqslant}

\def\GL{\operatorname{GL}}

\def\End{\operatorname{End}}

\def\dom{\operatorname{dom}}

\def\im{\operatorname{im}}

\def\Assoc{\operatorname{Assoc}}
\def\PBL{\operatorname{PBL}}
\def\codim{\operatorname{codim}}

\def\F{\mathbb{F}}

\def\cE{\EuScript{E}}

\def\0{\boldsymbol{0}}

\def\fra{\mathfrak a}

\begin{document}

\begin{center}
\Large\bf
One algebra of double cosets for a general linear group
over a finite field

\bigskip

\large\sc Yu. A. Neretin%
\footnote{Supported by the grant FWF (The Austrian Scientific Funds),
Project PAT5335224.}

\end{center}

\bigskip

{\small 
Let $\mathbb {F}_q$ be  finite field with $q$ elements. Let $\alpha\leqslant n$ be positive integers. Consider the general linear  group $\mathrm{GL}(\alpha+n, \mathbb {F}_q) $ and its subgroup $H(n)$, which fixes the first $\alpha$ basis elements in $\mathbb {F}_q^{\alpha+n}$. Denote $\mathcal{A}_n$ by the convolution algebra  of $H(n)$-biinvariant functions on $\mathrm{GL}(\alpha+n, \mathbb {F}_q) $. We describe  algebras $\mathcal{A}_n$ in terms of generators and relations and show that the family $\mathcal{A}_n$ admits a natural interpolation to  arbitrary complex $n$ (the field $\mathbb {F}_q$ and $\alpha$ are fixed).

}

\section{Introduction}

\COUNTERS

{\bf \punct Algebras of double cosets.%
\label{ss:algebra}}
Let $G$ be a finite group, let $\wh G$ be the set of its irreducible representations. For $\rho\in\wh G$,
we denote by $V_\rho$ a space of a representation $\rho$.
By  $\C[G]$ we denote the group algebra, its elements are
formal sums 
$$\sum_{g\in G}c_g g$$
 of elements of $G$ with complex coefficients.
Recall that
 $\C[G]$ is isomorphic to the direct sum
$$
\C[G]\simeq \bigoplus_{\rho\in \wh G} \End(V_\rho),
$$
where $\End(W)$ denotes the space of all linear operators in $W$.

Let $H$ be a subgroup in $G$. For any $\rho\in \wh G$, we denote by $V_\rho^H\subset V_\rho$
the subspace of  vectors fixed by all operators $\rho(h)$, where $h$ ranges in $H$
(this subspace can be trivial). Denote by $H\backslash G/H$ the space of double cosets
of $G$ with respect to $H$.

Denote by $\C[H\backslash G/H]\subset \C[G]$ the algebra of $H$-biinvariant elements of
the group algebra, i.e., sums $\sum c_g g$ 
such that for each $h\in H$ we have 
$$
h\cdot \Bigl( \sum_{g\in G} c_g g\Bigr)= \Bigl( \sum_{g\in G} c_g g\Bigr)
\cdot  h= \sum_{g\in G} c_g g.
$$  
Clearly, $\sum_{g\in G} c_g g\in \C[H\backslash G/H]$ if and only if the function $g\mapsto c_g$
is constant on double cosets $H\backslash G/H$.

This algebra also is a direct sum of matrix algebras,
$$
\C[H\backslash G/H]\simeq \bigoplus_{\rho\in \wh G} \End(V_\rho^H).
$$

Similarly, we can consider  
 a locally compact group $G$,  its compact subgroup $H$,
and the convolution algebra of $H$-biinvariant compactly supported complex-valued measures on $G$.

Various algebras of this type were topics of investigations after Gelfand's theorem about sphericity
\cite{Gel}, Mackey's works on induced representations \cite{Mack}, and Iwahori's constructions
of Hecke algebras and affine Hecke algebras \cite{Iwa1}-\cite{Iwa2}.
There is a collection of cases, when such algebras admit transparent descriptions, but usually they are 
heavy non-understandable objects. 

On the other hand, if $G$ is an infinite-dimensional
 group, and $H$ is a subgroup, then
quite often the set $H\setminus G/H$ has a semigroup structure. Such semigroups act
in spaces of $H$-fixed vectors in unitary representations. 
Moreover, such semigroups usually admit transparent
descriptions, see, e.g., \cite{Olsh-semi}, \cite{Ner-umn}.

 It appears that in the infinite-dimensional limit
a convolution in the group algebra quite often converges to a product of double cosets, see, \cite{Olsh-symm}, \cite{Olsh-semi},
\cite{Ner-concentration}.

The present paper is an attempt to extend a collection of hand
 algebras of double cosets. 

\sm 

{\bf \punct The purpose of the paper.} Let $\alpha\le n$ be positive integers. We consider the group
$G:=\GL(\alpha+n,\F_q)$ over a finite field $\F_q$ and its subgroup $H(n)$ consisting of transformations fixing
pointwise an $\alpha$-dimensional subspace in $\F_q^{\alpha+n}$. 

In this paper, we obtain  a description of the algebras
\begin{equation}
\C\bigl[H(n)\backslash \GL(\alpha+n,\F_q)/H(n)\bigr] 
\label{eq:double-1}
\end{equation}
in terms of generators and relations and construct an interpolating family of algebras depending on
a complex parameter  $n\in \C$.

Notice that the present work is parallel to the previous paper
\cite{Ner-PB} by the author devoted algebras
\begin{equation}
\C\bigl[S_n\backslash S_{\alpha+n}/S_n\bigr]
\label{eq:double-2}
\end{equation}
for symmetric groups, but description of the algebra \eqref{eq:double-2}
is more explicit than in our case.

Notice, that in both cases the limit theories as $n\to\infty$ are relatively simple. Namely, a limit pass for \eqref{eq:double-1}
leads to  the group of all permutations of a countable set (this is explained in \cite{Olsh-symm}
 and \eqref{eq:double-1} leads to the group of all linear operators  in a countable-dimensional linear space over $\F_q$
 (unitary representations of the latter group were classified in \cite{Tsa}). 

\section{Statements}

\COUNTERS

{\bf \punct Partial linear bijections.}
Let $\F_q$ be a finite field with $q$ elements.
 Below $q$ is fixed, so we write $\F$ instead of $\F_q$. 
Denote by $\GL(m,\F)$ the group  of invertible matrices over $\F$
of order $m$.  Number of its elements is
\begin{equation}
\# \GL(m,\F)=(q^m-1)(q^{m}-q)\dots(q^m-q^{m-1}).
\label{eq:GL}
\end{equation} 
We introduce the notation
\begin{equation}
[q;l]:=(q^l-1)(q^{l-1}-1)\dots (q-1),
\end{equation}
so 
\begin{equation}
\# \GL(m,\F)=[q;m]\cdot q^{m(m-1)/2}.
\end{equation}

\sm 

Let $V=\F^\alpha$ be a finite-dimensional linear space over $\F$.
A {\it linear partial bijection} $\lambda:V\to V$
is an invertible  linear operator
from a subspace $Y\subset V$ to a subspace $Z\subset V$.
We say that $Y$ is the {\it domain} $\dom\lambda$ of $\lambda$
 and $Z$ is the {\it image} $\im \lambda$ of $\lambda$.
 We say that the rank
 of $\lambda$ is 
 $$
 \rk \lambda:=\dim \dom \lambda=\dim \im\lambda.
 $$

We denote by $\PBL(\alpha,\F)$ the set of all partial linear bijection
$\F^\alpha\to\F^\alpha$.

For a subspace $S\subset \F^\alpha$  we denote by
$\cT[S]$ the partial linear bijection, whose domain is $S$
and which is identical on $S$.

\sm 

{\sc Remark.} Recall that partial linear bijections form a semigroup,
their multiplication is the natural multiplication of 
partially defined maps. The elements $\cT[S]$ are 
the idempotents in $\PBL(\alpha,\F)$. In this paper, we
 do not use
the multiplicative structure.
\hfill $\square$

\sm

{\bf \punct Description of double cosets.}
Let 
$$\boxed{\alpha\le n}$$ be  positive integers. Below $\alpha$ is fixed
and $n$ vary. 
We consider the group $\GL(\alpha+n,\F)$
and write its elements as
$$X=\begin{pmatrix}
a&b\\c&d
\end{pmatrix},
$$
subdividing blocks $a$, $b$, $c$, $d$
to smaller blocks if this is necessary.
{\it We apply matrices to vector-columns}.

We define the subgroup
$H(n)\subset \GL(\alpha+n,\F)$ consisting of block matrices
of size $\alpha+n$ having the form
$$
\begin{pmatrix}1&u\\0&v\end{pmatrix}
.$$ 
Clearly,
$$
\#H(n)=\bigl(\# \GL(n,\F)\bigr)\cdot q^{\alpha n}.
$$

For a matrix $g=\begin{pmatrix}A&B\\C&D\end{pmatrix}$
we define an element $\Pi(g)\in \PBL(\alpha,\F)$ as the map $x\mapsto  Ax$
restricted to the subspace $\ker C$.

\begin{proposition}
The map $\Pi$ is constant on double cosets $H(n)\cdot g\cdot H(n)$
and determines a bijection
$$
\Pi:
H(n)\backslash \GL(\alpha+n,\F)/H(n)\to \PBL(\alpha,\F).
$$ 
\end{proposition}

Denote by 
$$
\cA[\alpha;n]:=\C\bigl[H(n)\backslash \GL(\alpha+n,\F)/H(n)\bigr]
$$
 the subalgebra of
the group algebra $\C\bigl[\GL(\alpha+n;\F)\bigr]$ 
consisting of functions constant 
on double cosets $H(n)\backslash \GL(\alpha+n,\F)/H(n)$.

We also define the 'dual' map
$$
\Pi^\circ: \PBL(\alpha,\F)\to 
\cA[\alpha;n]
$$
by
$$
\Pi^\circ(\lambda)=\sum_{Q\in \GL(\alpha,\F):\, \Pi(Q)=\lambda} Q.
$$

For each subspace $S\subset \F^\alpha$, we define
the following element of $\cA[\alpha;n]$:
\begin{equation}
\theta(S):=\frac{1}{\# H(n)} \Pi^\circ\bigl(\cT[S]\bigr)
=  \frac{1}{\# H(n)}\sum_{Q:\,\Pi(Q)=\cT[S]} Q.
\label{eq:theta-S}
\end{equation}
 
 \sm

{\bf\punct Algebras $\bbA[\alpha;\nu]$.} We introduce some notation:

\sm

--- Denote by $\cP$  the set of all subspaces in $\F^\alpha$
of codimension 1.

\sm

--- For $L\in \cP$ denote by $\Xi(L)$
the set of all $g\in \GL(\alpha,\F)$ fixing
$L$ pointwise. Clearly,
$$
\# \Xi(L)=q^{\alpha-1}(q-1).
$$

--- 
 For two different subspaces $L$, $M\in \cP$
 we choose an element 
$\gamma_{L,M}\in \GL(\alpha,\F)$
such that, 
$$\gamma_{L,M}\,L=M,\quad \gamma_{L,M}\,M=L,\quad  \gamma_{L,M}^2=1,$$
 and
 $\gamma_{L,M}$  fixes all points of $L\cap M$ (below this choice formally
is present in  formulas but actually the expressions are independent 
 on this choice).  
 
 \sm

{\sc Generators.}
Fix $\nu\in \C$.
We consider an associative algebra $\bbA[\alpha;n]$  whose
generators are

\sm

--- $A(g)$, where $g$ ranges in $\GL(\alpha,\F)$,

\sm

--- $\Theta(L)$, where $L$ ranges in $\cP$.

\sm

These generators satisfy the following two groups
\eqref{eq:g1g2}--\eqref{eq:hL},
 \eqref{eq:Theta-square}--\eqref{eq:relation-last} of relations.

\sm

{\sc  The first  group of relations.}
\begin{align}
A(g_1)\, A(g_2)&=A(g_1 g_2);
\label{eq:g1g2}
\\
A(g)\, \Theta(L)\,A(g^{-1})&=\Theta(gL);
\label{eq:gL}
\\
A(h)\, \Theta(L)&= \Theta(L) \qquad \text{if $h\in \Xi(L)$}.
\label{eq:hL}
\end{align}

Clearly, $A(1)$ is the unit of 
the algebra $\bbA(\alpha;\nu)$. Also, all elements $A(g)$ are invertible,
$A(g)^{-1}=A(g^{-1})$.

\begin{corollary}
\label{cor:kgh}
{\rm a)} If  $h\in \Xi(L)$, then
\begin{equation}
 \Theta(L)\, A(h)= \Theta(L).
\end{equation}

\sm

{\rm b)}
 Let $g\in \GL(\alpha,\F)$,
$h\in \Xi(L)$, $k\in \Xi(gL)$.
Then
\begin{equation}
A(kgh)\,\Theta(L)= A(g)\,\Theta(L).
\label{eq:kgh}
\end{equation}	
\end{corollary}

\sm 

{\sc The second group of relations.}
For $L\in \cP$, we have:
\begin{equation}
\Theta(L)^2= (q^{\nu+1}-2q+1)\,q^{\alpha-1}
\Theta(L)+(q^\nu-1)q^{\alpha}\sum_{h\in \Xi(L)} A(h).
\label{eq:Theta-square} 
\end{equation}
Let  $L$, $M\in \cP$ be different. 
\begin{multline}
[\Theta(L),\Theta(M)]=(q-1)\,q^{\alpha-2}
\times \\ \times
\Bigl( \bigl[A(\gamma_{L,M}),\Theta(L)\bigr]
+\sum_{T\in \cP, T\supset L\cap M, T\ne L,M}\bigl(\,\bigl[A(\gamma_{T,L}),\Theta(L)\bigr]
-\bigl[A(\gamma_{T,M}),\Theta(M)\bigr]\,\bigr)
\Bigr),
\label{eq:Theta-Theta}
\end{multline}
where  $[\cdot,\cdot]$ denotes the commutator.

Finally, let $L\ne M$ be elements of $\cP$. Let 
$\eta\in \GL(\alpha,\F)$
fixes $L\cap M$ pointwise.
Then 
\begin{equation}
\bigl(A(\eta)-A(1)\bigr)
\cdot \Bigl(\Theta(L)\,\Theta(M)-
\sum_{T,S\in \cP:\, T,S\supset L\cap M,\, T\ne M, S\ne L}
A(\gamma_{T,S})\,\Theta(T))\Bigr)
=0.
\label{eq:relation-last}
 \end{equation}
  
 By Corollary \ref{cor:kgh}, the right-hand sides 
 of  \eqref{eq:Theta-Theta} and \eqref{eq:relation-last}
 do not depend on a choice of $\gamma_{T,S}$. 

\begin{theorem}
	\label{th:dim}
{\rm a)} For all $\nu\in \C$
$$
\dim\bbA[\alpha;\nu]=\# \PBL(\alpha,\F).
$$

{\rm b)} For   different $\nu$ algebras $\bbA[\alpha;\nu]$
are (noncanonically) isomorphic except a finite collection of values%
\footnote{Actually, coefficients in our  relations depend only on $q$ and $q^\nu$} of $q^\nu$.
\end{theorem}

\begin{proposition}
	\label{pr:M}
The linear functional $\text{\sf M}:\bbA[\alpha;\nu]\to \C$
defined on generators by 
$$
\text{\sf M}(A(g))=1,\qquad \text{\sf M}(\Theta(L))=
(q^\nu-1)\,(q^{\nu-1}-1)\,q^{2\alpha-1}
$$ 
is a homomorphism of algebras.
\end{proposition}

{\bf \punct The homomorphism 
$\boldsymbol{\Phi:\bbA[\alpha;n]\to \cA[\alpha;n]}$.}
For any $g\in \GL(\alpha;\F)$, we consider
the element
$$
a(g):= \frac{1}{\# H_n}\,\Pi^\circ(g)= \frac{1}{\# H_n} \sum_{R\in H(n)} g\cdot R\,\,\in \cA[\alpha;n].
$$
For any subspace $L\subset \F^\alpha$
of codimension 1 we consider an element
$$
\theta(L):= \frac{1}{\# H_n}\,\Pi^\circ\bigl(\cT[S]\bigr)
=  \frac{1}{\# H_n} \sum_{R\in \GL(\alpha+n,\F):\Pi(R)
=\cT[L]}
Q.
$$

\begin{theorem}
	\label{th:hom}
{\rm a)} The correspondence 
$$A(g)\mapsto a(g),\qquad \Theta(g)\mapsto \theta(g)$$
determines a homomorphism $\Phi:\bbA[\alpha;n]\to \cA[\alpha;n]$.

\sm

{\rm b)} For integer  $n\ge \alpha$ this homomorphism is an isomorphism.  
\end{theorem}

\begin{proposition}
	\label{pr:n-isomorphic}
	For all integer $n\ge \alpha$ the algebras 
	$\cA[\alpha;n]$ are (noncanonically) isomorphic.
\end{proposition}

{\bf \punct The filtration in $\bbA[\alpha;\nu]$ and PBW-theorem.}
We say that a {\it monomial expression} in $\bbA[\alpha;\nu]$ is an arbitrary product of generators
$A(h)$ and $\Theta(L)$. We say that a {\it $\Theta$-degree} of a monomial
expression is the number of factors $\Theta(L_j)$ in the product.
By the relations \eqref{eq:g1g2} and \eqref{eq:gL},
we can transform any monomial equation to an expression
of the form
$$
A(g)\,\Theta(L_1)\dots \Theta(L_p)
$$
of the same degree.

We say that a {\it polynomial expression}
is a sum of monomial expressions. Its  {\it $\Theta$-degree}
is the maximum  of  $\Theta$-degrees of summands.
We say that the {\it $\Theta$-degree} $\deg_\Theta(X)$ of an element
$X$
of $\bbA[\alpha;\nu]$ is the minimal $\Theta$-degree of polynomial
expressions representing $X$. Clearly,
$$
\deg_\Theta(X+Y)\le \max\bigl(\deg_\Theta(X), \deg_\Theta(Y)\bigr),
\qquad \deg_\Theta(XY)\le \deg_\Theta(X)+\deg_\Theta(Y).
$$

Denote by $\bbA[\alpha;\nu]_k\subset \bbA[\alpha;\nu]$
 the subspace of all elements of $\Theta$-degree
$\le k$,
$$
\bbA[\alpha;\nu]_0\subset \bbA[\alpha;\nu]_1\subset \bbA[\alpha;\nu]_2\subset\dots
$$
Clearly, $\bbA\bigl[\alpha;\nu]_0= \C[\GL(\alpha,\F)\bigr]$.

Denote by $\gr(\bbA)_k$ the quotient spaces
$$
\gr(\bbA)_k=\gr\bigl(\bbA[\alpha;\nu]\bigr):=
\bbA[\alpha;\nu]_k/\bbA[\alpha;\nu]_{k-1}.
$$

\begin{theorem}
	\label{th:filtration}
{\rm a)} $\deg_\Theta(X)\le\alpha$.

\sm

{\rm b)} $\dim \gr(\bbA)_k$ coincides with the number
of elements of $\PBL(\alpha,\F)$ of rank $\alpha-k$.

\sm

{\rm c)} Let $N\subset \F^\alpha$ be a subspace of codimension
$k$. Let
$$
N=\bigcap_{j=1}^k L_j,\qquad \text{where $L_j\in \cP$.}
$$ 
Then $\wh\Theta(N):=\Theta(L_1)\dots \Theta(L_k)$ as an element
of $\gr(\bbA)_k$ does not depend on a choice of $L_1$, \dots, $L_k$. 

\sm

{\rm d)} An element $A(g)\, \wh\Theta(N)$ as an element of $\gr(\bbA)_k$
depends only on $N$ and the restriction of $g$ to $N$,
i.e., on the partial linear bijection $\lambda$ of rank
$\alpha-k$, which is the restriction of $g$
to $N$.

\sm 

{\rm e)} Elements  $A(g)\, \wh\Theta(N)$, enumerated  by
partial linear bijections of rank $\alpha-k$,
 form a basis in $\gr(\bbA)_k$

\end{theorem}

\section{Proofs}

\COUNTERS

{\bf \punct Distributions of matrix elements of $\GL(m,\F)$.}
Recall  how  to count elements of a group $\GL(m,\F)$.
Let $\ell_1$, \dots, $\ell_m$ be rows of an invertible matrix
of order $m$, they form a basis of $\F^m$.
So, let us choose a basis. An element $\ell_1$ is an arbitrary
nonzero vector in $\F^m$, it can be selected in $q^m-1$ ways.
An element $\ell_2$ is contained in the complement
 to the line $\F\cdot\ell_1$.
It can be selected in $q^m-q$ ways. An element $\ell_3$ 
is contained in the complement to the plane $\F \ell_1\oplus \F\ell_2$.
It can be selected in $q^m-q^2$ ways, etc.

We  consider $\GL(m,\F)$ as a measure space
equipped with the uniform probabilistic measure.

\begin{lemma} 
\label{l:distribution-matrix}
{\rm a)} For $g\in \GL(m,\F)$, the matrix element $g_{11}$ is distributed
as follows:

\sm

--- $g_{11}$ is 0 with probability $(q^{m-1}-1)/(q^m-1)$;
 
\sm 
 
--- other values of $g_{11}$ are uniformly distributed,
each element has
 probability $q^{m-1}/(q^m-1)$.

\sm 

{\rm b)} For $g\in \GL(m,\F)$, a sub-row $(g_{11}, g_{12},\dots g_{1k})$
 is zero with probability
$(q^{m-k}-1)/(q^m-1)$.
Any nonzero vector $(g_{11}, g_{12},\dots g_{1k})$
 has probability 
$q^{m-k}/(q^m-1)$.
\end{lemma}

{\sc Proof.} a) Consider the subset in $\GL(m,\F)$, where $g_{11}=0$.
The vector $(g_{12}, \dots,g_{1m})$ is nonzero and can be selected in
$p^{m-1}-1$ ways. The remaining rows  $\ell_2$, \dots, $\ell_m$
 can be selected in 
 $$\prod\nolimits_{j=2}^m (q^m-q^j)$$
ways as above. 

If $g_{11}=a\ne 0$, when $(g_{12}, \dots,g_{1m})$ can be arbitrary.
It can be selected in $q^{m-1}$ way and we repeat the end of the previous paragraph.

\sm

b) Proof is similar. \hfill $\square$  

\sm
 
{\bf \punct Counting partial linear bijections.}
The group $\GL(\alpha,\F)\times \GL(\alpha,\F)$
acts on the set $\PBL(\alpha,\F)$ of all partial linear bijections  by left and right multiplications,
$$
(g,h): \lambda\mapsto g\,\lambda\, h^{-1}.
$$
Clearly, orbits are enumerated by a rank of $\lambda$.

\begin{lemma}
The number of elements of the set of all partial
linear bijections 
of rank $\alpha-\rho$
is
\begin{equation}
\sigma_\rho=
\frac{\bigl(\#\GL(\alpha,\F)\bigr)^2}
{\bigl(\#\GL(\rho,\F)\bigr)^2 \bigl(\#\GL(\alpha-\rho,\F)\bigr)
\cdot q^{2\rho(\alpha-\rho)}}.
\label{eq:number0}
\end{equation}
\end{lemma}

{\sc Proof.}
 We consider a coordinate subspace 
$\F^{\alpha-\rho} \subset \F^\alpha$
consisting of vectors (columns) $(x_1, \dots, x_{\alpha-r},0,\dots,0)^t$,
where $t$ denotes the transposition.
Consider the corresponding idempotent $\cT[\F^{\alpha-\rho}]$
in $\PBL(\alpha,\F)$. We must find its stabilizer in 
$\GL(\alpha,\F)\times \GL(\alpha,\F)$. 

 Let $g$, 
$h\in \GL(\alpha,\F)$. 
If $g\cdot \cT[\F^{\alpha-\rho}]\cdot h^{-1}=\cT[\F^{\alpha-\rho}]$,
then matrices $g$, $h$ have the following $((\alpha-\rho)+\rho)$-block structure
$$
g=\begin{pmatrix}
A&B\\0&C
\end{pmatrix}, \qquad
h=\begin{pmatrix}
A^{-1}&0\\P&Q
\end{pmatrix},
$$      
where
$$
A\in \GL(\alpha-\rho,\F),\quad C, Q\in \GL(\rho,\F),\quad
\text{$B$,  $P$, are arbitrary,}
$$
and this implies our statement.
\hfill $\square$

\sm

{\bf \punct Counting elements of double cosets.}

\begin{lemma}
Let a partial linear bijection has rank $\alpha-\rho$.
Then the corresponding double coset  has
\begin{equation}
\kappa_\rho=
\frac{\bigl(\# \GL(n,\F)\bigr)^2\,\cdot\, q^{2\alpha n}}
{\bigl(\#\GL(n-\rho;\F)\bigr)\,\cdot\,q^{(\alpha+\rho)(n-\rho)}}
\label{eq:number}
\end{equation}
elements.
\end{lemma}
In particular,
\begin{align}
\kappa_1&=(q^n-1)\,q^\alpha \cdot (\# H(n));
\label{eq:number1}
\\
\kappa_2&=(q^n-1)(q^{n-1}-1)\,q^{2\alpha+1} \cdot (\# H(n)).
\label{eq:number2}
\end{align}

{\sc Proof.}
As in the previous proof, we take the double coset
corresponding to the idempotent $\cT[\F^{\alpha-\rho}]$.
We take the following representative:
\begin{equation}
J_\rho=
\left(
\begin{array}{cc|cc}
1&0&0&0\\
0&0&1&0\\
\hline
0&1&0&0\\
0&0&0&1
\end{array}
\right),
\label{eq:J-rho}
\end{equation}
 sizes of blocks are 
\begin{equation} 
 \bigl((\alpha-\rho)+\rho\bigr)+\bigl(\rho+(n-\rho)\bigr)
 .
 \label{eq:size-rho}
 \end{equation}
 Clearly, $J_\rho^2=1$.
 The group $H(n)\times H(n)$ acts on $\GL(\alpha+n;\F)$ 
 by left and right multiplications,
 $$
(Q,R):\, X\mapsto QXR^{-1},\qquad X\in \GL(\alpha+n;\F),\,\, Q,R\in H(n), 
 $$
orbits are double cosets. Let us describe the stabilizer
of $J_\rho$.
Its elements satisfy 
$$
QJ_\rho R^{-1}=J_\rho,\quad\text{or, equivalently,}
\quad J_\rho  QJ_\rho=R. 
$$
Therefore, we must find number of 
all $Q\in H(n)$ such that $J_\rho  Q J_\rho\in H(n)$.
Representing $Q$ as a block matrix of size
\eqref{eq:size-rho}:
\begin{equation}
Q=\left(
\begin{array}{cc|cc}
1&0&b_{11}&b_{12}\\
0&1&b_{21}&b_{22}\\
\hline
0&0&d_{11}&d_{12}\\
0&0&d_{21}&d_{22}
\end{array}
\right),
\label{eq:q}
\end{equation}
we get
$$
J_\rho QJ_\rho=\left(
\begin{array}{cc|cc}
1&b_{11}&0&b_{12}\\
0&d_{11}&0&d_{12}\\
\hline
0&b_{21}&1&b_{22}\\
0&d_{21}&0&d_{22}
\end{array}
\right).
$$
Therefore, $b_{11}$, $b_{21}$, $d_{21}$ are zeros
and $d_{11}=1$. Hence,
\begin{equation}
Q=\left(
\begin{array}{cc|cc}
1&0&0&b_{12}\\
0&1&0&b_{22}\\
\hline
0&0&1&d_{12}\\
0&0&0&d_{22}
\end{array}
\right).
\label{eq:Q-1}
\end{equation}
The number of such $Q$ is the denominator
in \eqref{eq:number}, the numerator is
the number of elements in $H(n)\times H(n)$.
\hfill $\square$

\sm

{\sc Remark.} So, for $\sigma_\rho$ and $\kappa_\rho$ defined
by \eqref{eq:number0} and \eqref{eq:number}, we have
$$
\sum_{\rho=0}^\alpha \sigma_\rho \kappa_\rho=\#\GL(\alpha+n,\F).
$$
This is a special case of $q$-Chu--Vandermonde identity, see, e.g., 
\cite{AAR}, Ex.~10.4.b.
\hfill $\square$

\sm

{\bf \punct The verification of the relation \eqref{eq:hL}.%
\label{ss:eq-hL}}
The relations \eqref{eq:g1g2}--\eqref{eq:gL} for
$a(g)$ and $\theta(L)$ are obvious.

  To verify
 \eqref{eq:hL}, without loss of generality we can assume that
 $L$ is the coordinate subspace $L:=\F^{\alpha-1}\subset \F^\alpha$.
The double coset corresponding to $\cT[L]$ consists of invertible
block $\bigl((\alpha-1)+1+n\bigr)$-matrices of the form
\begin{equation} 
\left(\begin{array}{cc|c}
1&a_{12}& b_1\\
0&a_{22}&b_2\\
\hline
0&c_2& d
\end{array}
\right),\qquad
\text{where $c_2\ne 0.$}
\label{eq:??}
\end{equation}
 The group $\Xi(\F^{\alpha-1})\subset \GL(\alpha,\F)$ consists of
 matrices
 $$
 g=\left(\begin{array}{cc|c}
1&u_{12}& 0\\
0&u_{22}& 0\\
\hline
0&0& 1
\end{array}
\right),\qquad \text{where $u_{22}\ne 0.$}
 $$
Now let $R$ range in the set of matrices 
of the form \eqref{eq:??}. Then the map
$R\mapsto gR$ is a permutation of elements of our
double coset and does not change the sum
$$
\sum_{R:\,\Pi(R)=\cT[L]} R \qquad\in \cA[\alpha;n].
$$

\sm

{\bf \punct The proof of Corollary \ref{cor:kgh}.}
{\sc The statement a).}
$$
\Theta(L)\, A(h)=A(h)\, A(h^{-1})\, \Theta(L)\, A(h)=
A(h)\, \Theta(h^{-1} L)= A(h)\, \Theta(L)=\Theta(L).
$$

{\sc The statement b)}
Applying \eqref{eq:g1g2} and \eqref{eq:hL},
we get
$$
A(kgh)\,\Theta(L)=A(k)\,A(g)\,A(h)\,\Theta(L)=
A(k)\,A(g)\,\Theta(L).
$$
Next,
\begin{multline*}
A(k)\,A(g)\,\Theta(L)=A(k)\,\Bigl(A(g)\,\Theta(L)\, A(g^{-1})\Bigr)\, A(g)
=A(k)\,\Theta(gL)\, A(g)=
\\=\Theta(gL)\, A(g)
=A(g)\,\Theta(L)\, A(g^{-1})A(g)
=A(g)\,\Theta(L).
\end{multline*}

{\bf\punct The homomorphism 
$\boldsymbol{\text{\sf M}: \cA[\alpha;n]\to \C}$.}
For any finite group $G$ we have a canonical homomorphism
$$
\text{\sf M}: \C[G]\to\C
$$
defined by 
$$
\text{\sf M}\Bigl(\,\sum_{g\in G} c_g g\,\Bigr)=\sum_{g\in G} c_g.
$$
We can restrict it to any algebra $\C[H\backslash G/H]$.
In particular, we apply this remark to our case.
We have
$$
\text{\sf M}(A(g))=1,\qquad \text{\sf M}(\theta(L))=(q^n-1) q^{\alpha-1}.
$$

Clearly, for elements $\theta(S)$ defined by \eqref{eq:theta-S}
 we have
$$
\text{\sf M}(\theta(S))=\frac{\kappa_{\mathrm{codim}\, S}}{\# H(n)}.
$$

{\bf \punct The verification of the relation (\ref{eq:Theta-square}).}
We assume that  $L$ is  the coordinate subspace 
$\F^{\alpha-1}\subset \F^\alpha$
and take the matrix $J_1$, see \eqref{eq:J-rho}, as a representative
of the double coset corresponding $\cT[L]$.
To evaluate a convolution $\theta(L)\theta(L)$
we must find the distribution of double cosets
$$ H(n) J_1 Q J_1 H(n), \qquad \text{where $Q$ ranges in $H(n)$.}$$
For $Q$ of the form \eqref{eq:q} we have
 $$
 J_1 Q J_1=
\left( \begin{array}{cc|cc}
1&b_{11}&0&b_{12}\\
0&d_{11}&0&d_{12}\\
\hline
0&b_{21}&1&b_{22}\\
0&d_{21}&0&d_{22}
\end{array}
\right),\qquad \text{where $\det\begin{pmatrix}
d_{11}&d_{12}\\
d_{21}&d_{22}
\end{pmatrix}\ne 0$.}
$$

If $\begin{pmatrix}b_{21}\\d_{21}\end{pmatrix}\ne 0$, then
we fall again to $ H(n) J_1 Q J_1 H(n)$.  Otherwise 
we get the double coset
\begin{equation}
H(n)
\left( \begin{array}{cc|cc}
1&b_{11}&0&b_{12}\\
0&d_{11}&0&d_{12}\\
\hline
0&0&1&b_{22}\\
0&0&0&d_{22}
\end{array}
\right)
H(n)=
H(n)
\left( \begin{array}{cc|cc}
1&b_{11}&0&0\\
0&d_{11}&0&0\\
\hline
0&0&1&0\\
0&0&0&1
\end{array}
\right)
H(n).
\end{equation}
The submatrix 
$\begin{pmatrix}
1&b_{11}\\
0&d_{11}\\
\end{pmatrix}
$
ranges in $\Xi[L]$.
Therefore,
\begin{equation}
\theta(L)^2=s\cdot\theta(L)+\sum_{g\in \Xi(L)} t_g\, a(g)
\label{eq:theta-s-t}
\end{equation}
for some positive coefficients $s$, $t_g$.
Applying the homomorphism 
$$\text{\sf M}:\,\cA[\alpha;n]\to \C$$
to both sides, we get
\begin{equation}
(q^n-1)^2\, q^{2\alpha}=s\cdot (q^n-1)\, q^{\alpha}+ 
\sum_{g\in \Xi(L)} t_g.
\end{equation}

Let us evaluate the probability of the event
\begin{equation}
\begin{pmatrix}b_{21}\\d_{21}\end{pmatrix}=0.
\label{eq:event}
\end{equation} 
The element $b_{21}\in F$ is uniformly distributed,
so $b_{21}= 0$ with probability $1/q$.
The element $d_{21}$ is zero with probability $\frac{q-1}{q^n-1}$.
These two events are independent. So, we multiply these probabilities
and come to
$$\sum t_g= (q^n-1)^2\, q^{2\alpha} \cdot
 \Bigl(\frac 1q \cdot \frac{q-1}{q^n-1}\Bigr)
 =(q^n-1)\, q^{2\alpha-1}\, (q-1).
 $$
 Next, under the condition \eqref{eq:event} we have
$$
 Q=\left(
\begin{array}{cc|cc}
1&0&b_{11}&b_{12}\\
0&1&0&b_{22}\\
\hline
0&0&d_{11}&d_{12}\\
0&0&0&d_{22}
\end{array}
\right).
 $$
 Clearly, possible pairs $(b_{11}, d_{11})$ are uniformly
 distributed, $b_{11}\in \F^{\alpha-1}$, $d_{11}\in \F\setminus 0$.
  So all $t_g$ are equal
 and we get
 $$
 t_g=(q^n-1)q^{\alpha}.
 $$
Since we know $\sum t_g$ in \eqref{eq:theta-s-t},
we know also $s$.


\sm

{\bf\punct Verifying of the relation (\ref{eq:Theta-Theta}).}

\begin{proposition}
\label{pr:ThetaTheta}
Let $N\subset \F^\alpha$ be a subspace of codimension 2.
Denote
$$
$$
Then for different $L$, $M\subset \in\cP$ we have  
\begin{equation}
\theta(L)\,\theta (M)=
 \theta(L\cap M)+ q^{\alpha-2}(q-1)
\sum_{S,T\in \cP:\, S,T\supset N, S\ne M, T\ne L}
\beta_{S,T}\, \theta(S),
\label{eq:Theta-Theta-1}
\end{equation}
where elements $\beta_{S,T}\in \GL(\alpha,\F)$
are chosen from the conditions:

\sm

--- $\beta_{S,T}$ fixes $N$ pointwise;

\sm

--- $\beta_{S,T} S=T$.
\end{proposition}

{\sc Remark.} Elements $\beta_{S,T}$ can be chosen in many ways.
However, by Corollary \ref{cor:kgh} each summand
$\beta_{S,T} \Theta(S)$ does not depend on this choice.
\hfill $\boxtimes$ 

\sm

{\sc Proof.}
Without loss of the generality we can assume that
$L\subset \F^\alpha$ is the coordinate subspace 
$\F^{\alpha-1}$ consisting of vectors 
$(x_1,\dots,x_{\alpha-1},0)$
and $M$ is the coordinate subspace consisting of vectors 
$(y_1, \dots,y_{\alpha-2},0,y_\alpha)$.
We choose the following representatives of
these double cosets
\begin{equation}
I_L:=\left(\begin{array}{ccc|ccc}
1&0&0&0&0&0\\
0&1&0&0&0&0\\
0&0&0&0&1&0\\
\hline
0&0&0&1&0&0\\
0&0&1&0&0&0\\
0&0&0&0&0&1
\end{array}
\right),\qquad
I_M:=\left(\begin{array}{ccc|ccc}
1&0&0&0&0&0\\
0&0&0&1&0&0\\
0&0&1&0&0&0\\
\hline
0&1&0&0&0&0\\
0&0&0&0&1&0\\
0&0&0&0&0&1
\end{array}
\right),
\label{eq:ILIM}
\end{equation}
the size of matrices is $\bigl(((\alpha-2)+1+1)+(1+1+(n-2))\bigr)$.

We must examine double cosets 
$$
H(n)\cdot I_L Q I_M\cdot H(n),\qquad \text{where $Q$ ranges in $H(n)$.} 
$$
We write $I_L Q I_M$ explicitly and get
\begin{equation}
I_L
\left(
\begin{array}{ccc|ccc}
1&0&0&b_{11}&b_{12}&b_{13}\\
0&1&0&b_{21}&b_{22}&b_{23}\\
0&0&1&b_{31}&b_{32}&b_{33}\\
\hline
0&0&0&d_{11}&d_{12}&d_{13}\\
0&0&0&d_{21}&d_{22}&d_{23}\\
0&0&0&d_{31}&d_{32}&d_{33}
\end{array}
\right)
 I_M=
\left(
\begin{array}{ccc|ccc} 
1&b_{11}&0&0&b_{12}&b_{13}\\
0&b_{21}&0&1&b_{22}&b_{23}\\
0&d_{21}&0&0&d_{22}&d_{23}\\
\hline
0&d_{11}&0&0&d_{12}&d_{13}\\
0&b_{31}&1&0&b_{32}&b_{33}\\
0&d_{31}&0&0&d_{32}&d_{33}
\end{array}
\right),
\end{equation}
the sizes of matrices are $(\alpha-2)+1+1+1+1+(n-2))$. 
Let us examine  $\Pi(I_L Q I_M)\in\PBL(\alpha,\F)$.
We have
\begin{equation}
\dom \Pi(I_L Q I_M)=
\ker \begin{pmatrix}
0&d_{11}&0\\
0&b_{31}&1\\
0&d_{31}&0
\end{pmatrix}.
\label{eq:for-rank}
\end{equation}
Denote the last matrix by $\tau$. Its size is
 $\bigl((\alpha-2)+1+1\bigr)
\times\bigl(1+1+(n-2)\bigr)$, and therefore its rank $\le 2$,
the subspace $N=L\cap M$ is contained in the kernel.
If
$$
\begin{pmatrix}
d_{11}\\d_{31}
\end{pmatrix}\ne 0,
$$
then the rank is 2 and $\ker(\tau)=N$. Since $d\in\GL(\alpha+n,\F)$,
the probability of this event is
\begin{equation}
1-\frac{q-1}{q^n-1}=\frac{q(q^{n-1}-1)}{q^n-1}.
\label{eq:1-}
\end{equation}
If $d_{11}=0$, $d_{31}=0$, then we have the equation
$$
\begin{pmatrix}
0&0&0\\
0&b_{31}&1\\
0&0&0
\end{pmatrix}
\begin{pmatrix}
x_1\\x_2\\x_3
\end{pmatrix}=0,
$$
i.e., 
\begin{equation}
x_3=- b_{31} x_2, \qquad \text{$x_1$ is arbitrary.}
\end{equation}
So, this subspace, say $S$, has codimension 1, and $S\ne M$.
The coordinate $b_{31}\in\F$ on the set
\begin{equation} 
 \bigl\{Q\in H(n)\bigl|d_{11}=0, d_{31}=0\bigr\}
 \label{eq:curly} 
 \end{equation}
 is uniformly distributed, and therefore subspaces $S$ are uniformly
 distributed.
 
Next, we evaluate the map 
$\Pi(I_L Q I_M)$, 
$$
\begin{pmatrix}
1&b_{11}&0\\
0&b_{21}&0\\
0&d_{21}&0
\end{pmatrix}
\begin{pmatrix}
x_1\\
x_2\\
-b_{31} x_2
\end{pmatrix}
=\begin{pmatrix}
x_1+b_{11} x_2\\
b_{21} x_2\\
d_{21} x_2
\end{pmatrix}
.$$
On the set \eqref{eq:curly}
we have $d_{21}\ne 0$,
so $T:=\im \Pi(I_L Q I_M)\ne L$.
A subspace 
 $T$ is determined by the ratio $b_{21}/d_{21}$,
and $T\subset N$ can be an arbitrary  $(\alpha-1)$-dimensional
subspace except $L$. After the fixation of $b_{31}$ and
$b_{21}/d_{21}$ we can get an arbitrary linear map $S\to T$,
which is fixed on $N$, and such maps also are
uniformly distributed. 

So, we have
$$
\theta(L)\, \theta(M)= s\cdot \theta(L\cap M)+ t\cdot
\sum_{S,T\supset N: S\ne M, T\ne L}
\beta_{S,T} \theta(S),
$$
where $\beta_{S,T}$ are chosen as in the formulation
of the proposition. Applying the homomorphism $\text{\sf M}$
to both sides, we get
$$
\Bigl((q^n-1) q^\alpha\Bigr)^2=s\cdot (q^n-1)(q^{n-1}-1) q^{2\alpha+1}
+ (t q^2)(q^n-1) q^\alpha.
$$
Keeping in mind \eqref{eq:1-}, we observe
$$
\frac{q(q^{n-1}-1)}{q^n-1}
\Bigl((q^n-1) q^\alpha\Bigr)^2=s\cdot (q^n-1)(q^{n-1}-1) q^{2\alpha+1},
$$
and therefore $s=1$. Similarly,
$$
\frac{q-1}{q^n-1} \Bigl((q^n-1) q^\alpha\Bigr)^2
=t\cdot q^2\cdot(q^n-1) q^\alpha,
$$
and $t=(q-1)q^{\alpha-2}$.
\hfill $\square$

\sm

{\sc The proof of the relation (\ref{eq:Theta-Theta}).}
By Proposition \ref{pr:ThetaTheta} we have
\begin{multline}
\frac{q^{2-\alpha}}{q-1}
\bigl(\theta(L)\theta(M)-\theta(M)\theta(L)\bigr)
=\\=
\sum_{S,T\in \cP: S,T\supset N: S\ne M, T\ne L}
a(\gamma_{S,T}) \theta(S)- 
\sum_{S',T'\in \cP: S',T'\supset N: S'\ne L, T'\ne M}
a(\gamma_{S',T'}) \theta(T').
\label{eq:long}
\end{multline}
Next, we notice that for each pair $(S,T)$ we can choose $\gamma_{S,T}$
being an involution, $\gamma_{S,T}^2=1$, and we can assume that 
$\gamma_{S,T}=\gamma_{T,S}$.
Then for each $S$, $T\ne L$, $M$ both sums in the right hand side of 
\eqref{eq:long} contain a summand $\gamma_{T,S}\theta(S)$. Removing 
these summands, we see that the first sum in \eqref{eq:long}
reduces to
$$
a(\gamma_{M,L})\theta(L)+
\sum_{T\ne L,M} a(\gamma_{T,L})\theta(L) +
\sum_{S\ne L,M} a(\gamma_{S,M})\theta(S). 
$$
We can represent the last summand as $\sum_{T}\theta(M)A(\gamma_{T,M})$.

In the second sum in \eqref{eq:long} we get
$$
a(\gamma_{M,L})\theta(M)+
\sum_{T\ne L,M} a(\gamma_{T,M})\theta(M) +
\sum_{S\ne L,M} a(\gamma_{S,L})\theta(S). 
$$
The last summand is $\sum_S \theta_L \gamma_{S,L})$.
The difference gives our relation.
\hfill $\square$

\sm

Below in Subsect. \ref{ss:dim} we need the following lemma.

\begin{lemma}
\label{l:thetatheta}
Let $L$, $M$ be subspaces in $\F^\alpha$, $L+M=\F^\alpha$.
Then $\theta(L)\theta(M)$ has the form 
\begin{equation}
\theta(L) \theta(M)=\sigma \cdot \theta(L\cap M)+
\sum_{\lambda \in \PBL(\alpha,\F):\, \rk \lambda\ge \max(\dim L, \dim M)}
t_\lambda \Pi^\circ(\lambda)
\label{eq:thetathetatheta}
,\end{equation}
where $\sigma>0$, $t_\lambda\ge 0$.
\end{lemma}

{\sc Proof.} Without loss of generality we can assume that $L$, $M$ are
coordinate subspaces, $L$ consists of vectors
$(x_1, \dots, x_{\alpha-l},0,\dots,0)$
and $M$ of vectors
$(y_1,\dots,y_{\alpha-l-m}, 0,\dots,0,y_{\alpha-l},\dots,y_{\alpha})$.
Now we can repeat the beginning of the proof of Proposition
\ref{pr:ThetaTheta}. We take matrices $I_L$, $I_M$, see \eqref{eq:ILIM},
of size $\bigl(((\alpha-l-m)+m+l)+(l+m+(n-l-m))\bigr)$.
The rank of $\Pi(I_L Q I_M)$ is determined by the
rank of  matrix
$$
C:=\begin{pmatrix}
0&d_{11}&0\\
0&b_{31}&1\\
0&d_{31}&0
\end{pmatrix},
$$
see \eqref{eq:for-rank},
$$
\rk \Pi(I_L Q I_M)+\rk C=\alpha.
$$
Clearly, $l\le \rk C\le l+m$.

If $\rk C=l+m$, then $\ker C=L\cap M$, and 
$\Pi(I_L Q I_M)=L\cap M$. Clearly, this happens with nonzero probability.

On the other hand, $\rk \Pi(I_L Q I_M)=\rk \Pi(I_M Q I_L)$,
therefore $m\le \rk C$.
\hfill $\square$

\sm

{\bf \punct The verification of relations (\ref{eq:relation-last}).}
In the notation of Proposition \ref{pr:ThetaTheta},  we have
$$
A(\eta) \theta(L\cap M)= \theta(L\cap M). 
$$
the proof is as Subsect. \ref{ss:eq-hL},
we only replace block matrices of size 
$((\alpha-1)+1+n)$ by block matrices of size $((\alpha-2)+2+n)$.
Therefore,
$$
(A(\eta)-A(1)) \theta(L\cap M)=0
$$
We express $\theta(L\cap M)$ from \eqref{eq:Theta-Theta-1}
and get \eqref{eq:relation-last}.

\sm

{\bf\punct The verification Proposition \ref{pr:M}.}
Only a compatibility of $\text{\sf M}$ with \eqref{eq:Theta-square} requires a verification,
 i.e., we must check the identity 
\begin{equation*}
\M(	\Theta(L))^2= (q^{\nu+1}-2q+1)\,q^{\alpha-1}
	\M(\Theta(L))+\sum_{h\in \Xi(L)} \M(A(h)).
\end{equation*}
This is straightforward.

\sm 

{\bf \punct Proof of Theorem \ref{th:filtration}.c.}
Let $X$, $Y\in \bbA[\alpha;\nu]$ be contained in $\bbA[\alpha;\nu]_k$.
We write
$$
X\equiv_k Y
$$ 
if they are equal as elements of  $\gr(\bbA)_k$. By \eqref{eq:Theta-square} and \eqref{eq:Theta-Theta}, we have
\begin{equation}
\Theta(L)^2 \equiv_2 0,\qquad
\Theta(L)\Theta(M)\equiv_2 \Theta(M)\Theta(L).
\label{eq:l1l2}
\end{equation}

\begin{lemma}
	Let $L_1$, $L_2$, $M\in \cP$ be pairwise distinct, let $M\supset L_1\cap L_2$.
Then
$$
\Theta(L_1)\,\Theta(L_2)\equiv_2 \Theta(M)\,\Theta(L_2).
$$	
\end{lemma}

{\sc Proof.} The group $\Xi(L_2)$ acts transitively on the set of lines in $\F^\alpha$ 
transversal to $L_2$. Since each subspace $L_2$ and $M$ is generated by $L_1\cap L_2$
and such a line, there is $h\in \Xi(L_2)$ such that $hL_2=M$. Then
\begin{multline*}
\Theta(L_1)\Theta(L_2)\equiv_2 
\Theta(L_2)\Theta(L_1) \equiv_2 
A(h)\Theta(L_2)\Theta(L_2) \equiv_2 
A(h)\Theta(L_1)\Theta(L_2) \equiv_2 
\\ \equiv_2 
A(h)\Theta(L_1)\Theta(L_2)\equiv_2 \Bigl(A(h)\Theta(L_1)A(h)^{-1}\Bigl) \cdot \Bigl(A(h)\Theta(L_2)\Bigr)\equiv_2
\\ \equiv_2
\Theta(hL_2)\cdot \Theta(L_2)\equiv_2\Theta(M)\Theta(L_2),
\end{multline*}
here we use \eqref{eq:l1l2}, \eqref{eq:relation-last},  
\eqref{eq:gL},  and \eqref{eq:hL}.
	\hfill $\square$

\begin{corollary}
\label{cor:LLMM}
	Let $N\subset \F^\alpha$ be a subspace of codimension 2. Let $L_1$, $L_2$, $M_1$, $M_2\in\cP$
	and $L_1\cap L_2=N$, $M_1\cap M_2=N$. Then
	$$
	\Theta(L_1)\,\Theta(L_2)\equiv_2 \Theta(M_1)\,\Theta(M_2).
	$$
\end{corollary}
	
	{\sc Proof.} If $L_1$, $L_2$, $M_1$, $M_2$ are not pairwise distinct, the statement immediately follow from 
	the lemma. If they are distinct, then we apply lemma two times,
$$
\qquad
\Theta(L_1)\Theta(L_2)\equiv_2 \Theta(M_2)\Theta(L_2)\equiv_2 \Theta(L_2)\Theta(M_2)\equiv_2
\Theta(M_1)\Theta(M_2).
\qquad\square
$$



\sm

{\sc Proof of Theorem \ref{th:filtration}.c.}
So, let $L_1$, \dots, $L_k$, $M_1$, \dots, $M_k\in \cP$, and
$$N= L_1\cap\dots \cap L_k=M_1\cap\dots \cap M_k $$
 has codimension $k$ in $\F^\alpha$. We must show that
 \begin{equation}
 \Theta(L_1)\dots \Theta(L_k)\equiv_k \Theta(M_1)\dots \Theta(M_k).
 \label{eq:many-theta}
 \end{equation}

Let $\ell_j$, $m_j$  be linear functional on $\F^\alpha$
such that $\ker \ell_j=L_j$, $\ker m_j=M_j$
(so, we have two bases in the space $(\F^\alpha/N)'$
 dual to our $\F^\alpha/N$). Let 
$
\begin{pmatrix}
 u&v\\w&t
\end{pmatrix}\in \GL(2,\F).
$
Let us change the pair $L_i$, $L_j$ by 
$$
L_i'=\ker (u \ell_i+v \ell_j),\qquad
L_j'=\ker (w \ell_i+t \ell_j).
$$
By Corollary \ref{cor:LLMM}, such transformation does not change
the product of $\Theta(L_p)$ as an element of $\gr(\bbA)_k$.
Let $S_1$, $S_2$, \dots be linear transformations of $(\F^\alpha)'$
  of this type.
We can apply $S_1$, after this we can apply $S_2$ to a new basis, i.e.,
$S_1 S_2 S_1^{-1}$ in the old coordinates, and we get the transformation
$S_2S_1$. But clearly, such matrices generate the whole group 
$\GL(\alpha,\F)$. So we can transfer any basis in
the dual space to any basis. This proves \eqref{eq:many-theta}.

\begin{lemma}
Let $N$, $L_j$ be the same and $K\in \cP$ contains $N$.
Then
$$
\Theta(K) \Theta(L_1)\dots \Theta(L_k)\equiv_{k+1}0. 
$$ 
\end{lemma}

{\sc Proof.}
Indeed, by \eqref{eq:many-theta} without loss of generality
we can assume $L_1=N$, and we apply the relation \eqref{eq:Theta-square}.
\hfill $\square$

\sm

{\bf \punct Proof of Theorem \ref{th:filtration}.d.}

\begin{proposition}
Let $N$ and $L_1$, \dots, $L_k$ be the same as above. 
Let $\xi\in \GL(\alpha,\F)$ fix $N$ pointwise.
Then
$$
A(\xi)\, \wh\Theta(N)\equiv_k \wh\Theta(N).
$$
\end{proposition}

{\sc Proof.}  Without loss of generality we can assume that
that $L_k$, $L_{k-1}$, $L_{k-2}$, \dots are coordinate subspaces, consisting of vectors
$$
(x_1,\dots,x_{\alpha-1},0), \,\, (x_1,\dots,x_{\alpha-2},0, x_\alpha),\,\, (x_1,\dots,x_{\alpha-3},0, x_{\alpha-1}, x_\alpha), \dots
$$
respectively. The subgroup $Q_k$ in $\GL(\alpha,\F)$ fixing $N$ pointwise
consists of block matrices of size $((\alpha-k)+k)$ of the form
$$
\left(\begin{array}{c:c}
	1& u\\
	\hdashline
	0&v
\end{array}
\right).
$$
It is a semidirect product of $\GL(k,\F)$ (corresponding to the block $v$)
 and the Abelian subgroup $\F^{k(\alpha-k)}$ (corresponding to matrices $v$).

 A group $\Xi[L_j]\subset Q_k$ consists of block
$\bigl((\alpha-k)+(j-1)+1+(k-j)\bigr)$-matrices
of  the following structure
$$ 
\left(
 \begin{array}{c:ccc}
 	1&0&u_{k}&0\\
 	\hdashline
 	0&1&v_{1k}&0\\
 	0&0&v_{2k}&0\\
 	0&0&v_{3k}&1
 \end{array}
 \right),\qquad \text{where $v_{1k}\ne0$.}
$$
Obviously, the subgroups $\Xi[L_1]$, \dots, $\Xi[L_k]$,
generate $Q_k$. Indeed, these subgroups contain all Chevalley generators
for $\GL(\alpha-k,\F)\subset Q_k$. On the other hand
subgroups in $\Xi[L_j]$ consisting of matrices of the form 
$$ 
\left(
 \begin{array}{c:ccc}
 	1&0& *&0\\
 	\hdashline
 	0&1&0&0\\
 	0&0&1&0\\
 	0&0&0&1
 \end{array}
 \right)
$$
generate $\F^{k(\alpha-k)}$.

 Let $\eta\in\Xi[L_j]$
Then  we  can reoder products in $\gr(\bbA)_k$.
 \begin{multline*}
 A(\eta) \Theta(L_1)\dots \Theta(L_k)\equiv_k 
  A(\eta) \Theta(L_j)\Theta(L_1) \Theta(L_2)\dots \equiv_k
  \\ \equiv_k
  \Theta(L_j)\Theta(L_1) \Theta(L_2)\dots \equiv_k
\Theta(L_1)\dots \Theta(L_k).  
\qquad\qquad \square
 \end{multline*}
 
 \sm
 This statement implies Theorem \ref{th:filtration}.d.
 
 \sm 
 
 {\bf \punct The dimension of $\bbA[\alpha;n]$.%
 \label{ss:dim}}
 
 \begin{corollary}
 	\begin{equation}
 	\dim \bbA[\alpha;\nu]\le \# \PBL(\alpha,\F).
 	\label{eq:dim-dim}
 	\end{equation}
 \end{corollary}
 
 {\sc Proof.} Indeed, 
 $$
 \dim \bbA[\alpha;\nu]=\sum_{k=0}^\alpha \gr(\bbA)_k.
 $$
 On the other hand, each space $\gr(\bbA)_k$ contains a generating family of vectors enumerated
 by partial linear bijections of rank $\alpha-k$.
 \hfill $\square$

\begin{lemma}
Let $n\ge \alpha$. Then
$$
\dim \bbA[\alpha;n]= \# \PBL(\alpha,\F).
$$
\end{lemma} 

\begin{lemma} 
The homomorphism $\Phi:\bbA[\alpha;n]\to \cA[\alpha;n]$ is surjective. 
\end{lemma} 
 
 {\sc Proof.}
We prove this in the following form:

\sm

---  {\it for each $k$ the image of $\Phi$ contains all elements 
$\theta(L)$ with $\codim L\le k$.}

\sm

We prove this statement by induction in $k$. Assume that this valid for 
$k=m$. Let $\codim S=m-1$. We represent $S$ in the form $S=L\cap M$,
where $\codim L=m$, $\codim M=1$. We apply Lemma \ref{l:thetatheta}.
By the induction hypothesis, $\theta(L)$, $\theta(M)$ are contained
in $\im \Phi$. Therefore, $\theta(L)\theta(M)\in \im\Phi$.
By induction hypothesis, the term $\sum(\dots)$ in the right hand 
side of \eqref{eq:thetathetatheta} is contained in $\im \Phi$.  
Hence, $\theta(L\cap M)=\theta(S)\in \im\Phi$.
\hfill $\square$ 
 
\sm

 {\bf \punct The dimension of $\bbA[\alpha;\nu]$, where $\nu\in\C$.%
 \label{ss:dim1}}
For each $k$ and each subspace $N \subset \F^\alpha$ 
of codimension $k$ we fix
 an ordered collection $L_1^\circ$, \dots, $L_k^\circ\in \cP$ such that $N=\cap L_j$.
Also, for each  $\Lambda\in \PBL(\alpha,\F) $ we fix an element $g^\circ\in \GL(\alpha)$,
which coincides with $\lambda$ on $\dom \lambda$.
So, for any $\Lambda$ of rank $\alpha-k$ we get the following 
 element of $\bbA[\alpha;\nu]$
 \begin{equation}
 A(g^\circ)\Theta(L_1^\circ)\dots \Theta(L_k^\circ).
 \label{eq:basis}
 \end{equation}
 Clearly, for each $k$ images of such elements in $\gr(\bbA)_k$
 do not depend on the choice of $L_j^\circ$, $g^\circ$. If $\nu=n$, these images 
 form a basis in  $\gr(\bbA[\alpha;n])_k$. Therefore, the collection
\eqref{eq:basis} forms a basis  in $\cA[\alpha;n]$.
 
 To simplify notation, we denote elements of this basis by $e_\mu$.
 Products of basis elements have the form 
 $$
 e_\mu e_\nu=\sum \fra_{\mu,\nu}^\kappa (q;n)\, e_\kappa,
 $$
 where {\it $\fra_{\mu,\nu}^\kappa (q;n)$ are polynomial expressions in $q$, $q^{\alpha-2}$, and $q^n$}.
 In our considerations, $q$ and $q^{\alpha-2}$ are constants, and $q^n$ vary. So, we can regard 
 $\fra_{\mu,\nu}^\kappa$ as polynomials $\fra_{\mu,\nu}^\kappa(q^n)$ in an variable $q^n$ with constant coefficients.
 
 The associativity 
 $$ (e_\mu e_\nu) e_\kappa= e_\mu (e_\nu e_\kappa)$$
 of $\cA[\alpha;n]$ is equivalent to a collection of quadratic equations for coefficients 
 ({\it structure constants})
 (see, e.g. \cite{Maz})
 $\fra_{\mu,\nu}^\kappa(q^n)$,
 \begin{equation}
 \sum_\kappa \fra_{\mu,\nu}^\kappa(q^n) \,\fra_{\kappa,\pi}^\psi(q^n)=
 \sum_\kappa \fra_{\mu,\kappa}^\psi(q^n) \,\fra_{\nu,\pi}^\kappa(q^n),
  \label{eq:structure}
 \end{equation}
 equations are enumerated by collections $\kappa$, $\mu$, $\nu$, $\psi$.
 In both sides of each equation we have polynomials in $q^n$. We have an identity in countable number of points,
  $q^n=q^\alpha$, $q^{\alpha+1}$, $q^{\alpha+2}$, \dots. Therefore, this holds for all complex values of $q^n$.
  
  Thus we get an equidimensional family of associative algebras, 
  say $\cA[\alpha;\nu]$, holomorphically depending on
  the  complex parameter $\nu$ with structure constant polynomial in $q^\nu$.
  Each algebra  is generated by $a(g)$, $\theta(L)$, which are elements 
  of the basis \eqref{eq:basis}. The relations
  \eqref{eq:g1g2}--\eqref{eq:relation-last} are fulfilled for
   integer $\nu\ge\alpha$, and therefore they are
  fulfilled for al $\nu\in\C$.
  
  Therefore, $\cA[\alpha;n]$ is a quotient algebra of $\bbA[\alpha;n]$. But  inequality \eqref{eq:dim-dim}
  shows that the kernel is trivial. Therefore, 
  $$
  \bbA[\alpha;\nu]=\cA[\alpha;\nu].
  $$
  
  So, we proved Theorem \ref{th:dim}.a, Theorem \ref{th:hom}.b, Theorem \ref{th:filtration}.b,e.

\sm

{\bf \punct Proof of Theorem \ref{th:dim}.b. 
The curve $\boldsymbol{\nu\mapsto\bbA[\alpha;\nu]}$ in the space of algebras.}
Fix a space $\C^N$.
The space $\Assoc(\C^N)$ of structures of associative algebras on this space is an algebraic variety
defined by a system of quadratic equations \eqref{eq:structure}. The group $\GL(N,\C)$ acts
on this space by changing of bases.
Classes of isomorphisms if $N$-dimensional algebras correspond to orbits of $\GL(N,\C)$.
 The variety $\Assoc(\C^N)$ has many irreducible components,
  a few is known on this topic (see, e.g., \cite{Maz}, \cite{Ner-Lie}). 
  
  It is well-known
  and semi-obvious that 
  each semisimple algebra%
  \footnote{
  Recall that  a finite-dimensional semisimple associative algebra is a direct sum of matrix algebras.} $B$ determines a component of $\Assoc(\C^N)$, see. e.g., \cite{Fla}. Namely, there is a component 
  $\ov {\cO(B)}$,
  where $\GL(N,\C)$ has an open dense orbit $\cO(B)$, and this orbit corresponds to algebras isomorphic to $B$. The algebraic  variety 
  $\ov {\cO(B)}\setminus \cO(B)$ has codimension $\ge 1$
  in $\ov {\cO(B)}$. It can not contain other semisimple algebra $C$, since
  $\cO(C)$ is dense  in its own component $\ov{\cO(C)}$.

  Next, for any finite group $G$ and any subgroup $K$ the algebra 
  $\C[K\backslash G/K]$ is semisimple. 
  The curve $\nu\mapsto \bbA[\alpha;\nu]$ has a polynomial parametrization.
   Therefore it is algebraic and is completely contained in the component $\ov {\cO\bigl(\bbA[\alpha;n]\bigr)}$
   for each $n\ge \alpha$. In particular, all components 
   $\ov {\cO(\bbA[\alpha;n)]}$ coincide, and all algebras
 $\bbA[\alpha;n]$ are isomorphic. So, we get Proposition \ref{pr:n-isomorphic}.
  
  On the other hand, this curve  can have  only a finite number of points of intersection with 
  the  stratum $\ov {\cO(\bbA[\alpha;n])}\setminus\cO(\bbA[\alpha;n)])$. So, we get Theorem \ref{th:dim}.b.
  
\sm

We also get the following corollary. In the notation   of Subsect. \ref{ss:algebra} for each $n$
consider the collection
\begin{equation*}
	\cE_n:=
\bigl\{\dim V_\rho ^{H(n)}\bigr\},\qquad \text{where $\rho$ ranges in $\wh{\GL(\alpha+n,\F)\vphantom{\Bigl|}}$}.
\end{equation*}

\begin{proposition}
For all $n\ge\alpha$ collections of non-zero elements of $\cE_n$ is the same.
\end{proposition}

\noindent
\tt
University of Graz,
Department of Mathematics and Scientific computing;
\\
High School of Modern Mathematics MIPT;
\\
Moscow State University, MechMath. Dept;
\\
University of Vienna, Faculty of Mathematics.
\\
e-mail:yurii.neretin(dog)univie.ac.at
\\
URL: https://www.mat.univie.ac.at/$\sim$neretin/
\\
\phantom{URL:}
https://imsc.uni-graz.at/neretin/index.html

\end{document}